\documentclass[10pt]{article}
\usepackage{latexsym}
\usepackage{amsmath}
\usepackage{amssymb}
\usepackage{amsthm}

\newtheorem{theorem}{Theorem}
\newtheorem{proposition}[theorem]{Proposition}

\newtheorem{corollary}[theorem]{Corollary}

\newcommand{\R}{\mathbb{R}}
\newcommand{\N}{\mathbb{N}}
\newcommand{\Z}{\mathbb{Z}}
\newcommand{\Q}{\mathbb{Q}}
\newcommand{\C}{\mathbb{C}}
\newcommand{\e}{\mathrm{e}}
\newcommand{\im}{\mathrm{i}} 
\renewcommand{\O}{\mathrm{O}}

\title{Non-Holonomicity of Sequences Defined via Elementary Functions}
\author{Jason P. Bell \and Stefan Gerhold \and Martin Klazar \and Florian Luca}

\author{Jason P. Bell\footnote{Simon Fraser University, Vancouver. {\tt jpb at sfu.ca}},
Stefan Gerhold\footnote{Vienna University of Technology. {\tt sgerhold at fam.tuwien.ac.at}},
Martin Klazar\footnote{Charles University, Prague.  {\tt klazar at kam.mff.cuni.cz}},
Florian Luca\footnote{Instituto de Matem{\'a}ticas UNAM, Campus Morelia. {\tt fluca at matmor.unam.mx}}
}

\date{\today}

\begin{document}

\maketitle
\begin{abstract}
We present a new method for proving non-holonomicity of sequences, which is based
on results about the number of zeros of elementary and of analytic functions.
Our approach is applicable to sequences that are defined as the values of an elementary function at
positive integral arguments. We generalize several recent results; e.g., non-holonomicity of the logarithmic sequence
is extended to rational functions involving $\log n$.
Moreover, we show that the sequence that arises from evaluating the Riemann zeta function
at odd integers is not holonomic.
\end{abstract}

\begin{center}
\end{center}

\section{Introduction}

A sequence of complex numbers $(a_n)_{n\ge 1}$ is called {\em holonomic} (P-recursive) if there 
exist polynomials $p_k(x)\in\C[x]$, $k=0,1,\dots,d$, which are not all zero, such that
\begin{equation}\label{Prec}
p_0(n)a_n + p_1(n)a_{n+1}+ \dots +p_d(n)a_{n+d} = 0
\end{equation}
holds for every integer $n\ge 1$ (equivalently, for every integer $n\ge n_0$).
Their continuous counterpart are holonomic ($D$-finite) functions $g(z)$, which satisfy linear differential equations
\[
  p_0(z)g(z) + p_1(z)g'(z) + \dots + p_d(z)g^{(d)}(z) = 0
\]
with polynomial coefficients. For more information on holonomic sequences and functions,
especially in combinatorial enumeration, see Stanley~\cite{Stanley98}.
For the sake of clarity, in this paper we will use the word holonomic for sequences
and the word $D$-finite for functions.
Sequences that satisfy~\eqref{Prec} with $d=1$ are called hypergeometric.
An important property for our purposes is that the sum and the product of two holonomic sequences are holonomic.

The role of the holonomic sequences in the set of all complex sequences is reminiscent of the role
of the algebraic numbers in the set of all complex numbers.
Just like complex numbers are usually not algebraic, unless they are by design, a sequence
that is not obviously holonomic is usually not holonomic.
In both situations it can be a challenge, though, to come up with a proof.
Flajolet, Gerhold, and Salvy~\cite{FlGeSa:05} give an exhaustive survey of known non-holonomicity results.
In the present article we shall be interested in the following

\bigskip\noindent
{\bf Question.} If $f:\;\Omega\to\C$ is a ``nice'' function that is defined on a domain $\Omega\subset\C$ 
containing all natural 
numbers $\N=\{1,2,3,\dots\}$, when is it the case that the sequence of values of $f$ on $\N$, $(f(n))_{n\ge 1}$, is 
non-holonomic?

\bigskip\noindent
We seek general results and criteria yielding non-holonomicity of $(f(n))_{n\ge 1}$ for $f$ from large classes 
of functions ${\cal F}$, generalizing the results from Gerhold~\cite{Ge04} and some of the results of Flajolet, Gerhold, and 
Salvy~\cite{FlGeSa:05}. In the first paper, $(\log n)_{n\ge 1}$ was proved conditionally to be non-holonomic, and 
in the second paper an unconditional proof was given by an asymptotic machinery.
A simple proof was given also by Klazar~\cite{Kl05}. Here, we generalize the method from the latter paper
and arrive at non-holonomicity results for many sequences for which closed-form expressions are available.

In the following section, we explain our proving method. Section~\ref{se:res} is devoted
to straightforward applications of our approach. In Section~\ref{se:exp} we develop
our arguments further in several directions to establish the non-holonomicity of some sequences
involving the exponential and factorial functions.
We proceed by presenting several unrelated results, e.g. about algebraic sequences and interlacement sequences,
in Section~\ref{se:misc}. 
Finally, we show in Section~\ref{se:zeta} that
the sequence arising from evaluating the Riemann zeta function at odd integers is not holonomic.
The method we use there is not related to our main approach.

\section{The Proving Method}\label{se:methods}

\bigskip
Before presenting our approach we briefly discuss the method of Flajolet, Gerhold, and Salvy~\cite{FlGeSa:05}.
They use the equivalence saying that $(a_n)_{n\ge 1}$ is holonomic iff the ordinary generating function 
$g(z)=\sum_n a_nz^n$ is $D$-finite~\cite{Stanley98}. The asymptotic behavior of a holonomic function near a singularity is constrained by a structure theorem~\cite[Theorem 2]{FlGeSa:05}, and Abelian theorems transfer the asymptotic behavior 
of $a_n$ as $n\to\infty$ to the asymptotic behavior of $g(z)$ near the singularity~\cite[Theorem 3]{FlGeSa:05}. 
Thus, if $(a_n)_{n\ge 1}$ (or any transform of it obtained by holonomicity-preserving tranformations) 
has an asymptotic behaviour that is transferred to ``forbidden'' asymptotics of the generating function,
then $(a_n)_{n\ge 1}$ is not holonomic. While requiring an initial asymptotic analysis of 
$(a_n)_{n\ge 1}$, which is a drawback, the real strength of this method is that this asymptotic behaviour is 
fully sufficient information, and a closed-form representation $a_n=f(n)$ is not needed to prove 
non-holonomicity. This has been demonstrated, e.g., on the sequence of primes~\cite{FlGeSa:05}. 
In contrast, our method relies heavily on the explicit representation $a_n=f(n)$.

We associate with a function $f(z)$ and a $(d+1)$-tuple of complex polynomials 
$p_0(z),\dots,p_d(z)$ the function
$$
F(z)=F(z;f,p_0,\dots,p_d):=\sum_{i=0}^d p_i(z)f(z+i).
$$
If $(f(n))_{n\ge 1}$ is holonomic, then, for some not-all-zero polynomials $p_0,\dots,p_d$,
the function $F(z)$ vanishes at 
$z=1,2,\dots$. 
If $f$ is real-valued in the positive reals (or at least in $\N$), then $(f(n))_{n\ge 1}$ is a 
sequence of real numbers, and by an algebraic argument~\cite{Ge04,Lipshitz89} 
the polynomials in~\eqref{Prec} can be assumed without loss of generality to have real coefficients.
This is important for the formulation of the following two properties for
classes of complex functions ${\cal F}$: If all $f\in\mathcal{F}$ are real-valued for real arguments,
we assume that the polynomials $p_0(z),\dots,p_d(z)$ have real coefficients.
\begin{description}
\item[Property (A)] If $f\in{\cal F}$ is not identically zero and not all (real or complex, see above) polynomials $p_0(z),\dots,p_d(z)$ are zero, then the function $F(z;f,p_0,\dots,p_d)$ is not identically zero.
\item[Property (B)] If $F(z)=F(z;f,p_0,\dots,p_d)$, $f\in{\cal F}$, vanishes at all $z\in\N$,
then $F(z)$ is identically zero.
\end{description}
It follows that any non-identically-zero function $f$ from a class ${\cal F}$ with properties~(A) and (B)
produces a non-holonomic sequence $(f(n))_{n\ge 1}$.
(For some of the classes $\mathcal{F}$ we are going to discuss, $1$ has to be replaced
by some other positive integer $n_0$.)

\bigskip
Let us have a look at some (classes 
of) functions having property~(A) and then at functions having property~(B).
Ultimately we want to have large classes of functions having simultaneously properties~(A) and (B).

\bigskip

{\bf Condition (A1). Meromorphic functions and singularities.} The class of functions $f$ which are meromorphic for $\Re(z)>0$ but have no extension that is meromorphic at $z=0$ 
has property (A). Indeed, if $p_0$ is a nonzero polynomial and 
$F(z;f,p_0,\dots,p_d)\equiv 0$, then
\begin{equation}\label{leftrightside}
f(z) = -p_0(z)^{-1} \sum_{i=1}^d p_i(z)f(z+i)
\end{equation}
for $\Re(z)>0$, which is impossible, because the right hand side is meromorphic at $z=0$ while the left hand side is not. 
This criterion can be generalized in an obvious way by replacing $0$ with any other point $z_0\in\C$. 
Examples of such functions are $\log z$, $\log\log z$, $z^{\alpha}$ with $\alpha\in\C\backslash\Z$, $\exp(\sqrt{z})$, 
$\exp(\sqrt{\log z})$, $\sin(1/z)$, $z^z$, $\arctan z$,  etc.
Condition~(A1) will be our main tool for establishing property~(A).
Note how the assumption that our sequences $a_n$ are defined via functions $f(z)$ enables us to work with
the asymptotic behaviour of $f(z)$ at points other than infinity.

\bigskip
{\bf Condition (A2). Logarithmic derivatives and shifts.}
Here we assume that $f'(z)/f(z)\in\C(z)$. Notice that this implies that $f$ is of the form 
\[
  f(z)=\exp(r(z))\prod_{i=1}^m (z-\alpha_i)^{c_i},
\] 
where $r(z)$ is a rational function, and $\alpha_i,c_i\in\C$.

\begin{proposition}\label{prop:jason}
The class of functions $f(z)$ that are analytic in a neighborhood of ${[1,{\infty}[}$
and satisfy $f'(z)/f(z)\in\C(z)$ but $f(z+k)/f(z)\not\in\C(z)$, $k\in\N$, has 
property (A).
\end{proposition}
\begin{proof}
Suppose, on the contrary, that we have a nontrivial relation $F(z;f,p_0,\dots,p_d)\equiv 0$,
with minimal $d$. Upon differentiating, this implies
\[
  0 = F'(z) = \sum_{i=0}^d \left(p_i'(z) + p_i(z)\frac{f'(z+i)}{f(z+i)} \right) f(z+i).
\]
We now eliminate the summand $i=0$ from the equations $F(z)=0$ and $F'(z)=0$.
The coefficients $b_i(z)$ in
\begin{equation*}
  \sum_{i=1}^d b_i(z)f(z+i) := F'(z)p_0(z) - F(z)\left( p_0'(z) + p_0(z) \frac{f'(z)}{f(z)} \right) \in \C(z)
\end{equation*}
are in $\C(z)$, hence they vanish by the minimality of $d$. For $i=d$ this reads as
\[
  0=p_0(z)\left(p_d'(z)+p_d(z)\frac{f'(z+d)}{f(z+d)}\right)
   - \left(p_0'(z)+p_0(z)\frac{f'(z)}{f(z)}\right)p_d(z).
\]
Dividing by $p_0(z)p_d(z)$ and integrating, we obtain
\[
  f(z+d)/f(z) = C p_0(z)/p_d(z) \in \C(z)
\]
for some constant $C$, a contradiction. Note that we have $p_0(z)p_d(z)\neq0$ by the minimality of $d$.
\end{proof}
Notice that the field $\C(z)$ in the proposition can be replaced by an arbitrary field $\mathbb{K}\supset\C(z)$ of functions meromorphic at infinity, provided that it is closed under differentiation and shift.
For instance, this generalization can be applied to $f(z)=\exp\exp(1/z)$ with
$\mathbb{K}=\C(z, \e^{1/z}, \e^{1/(z+1)},\dots)$.

We remark that we do not know of examples where condition~(A2) works and condition~(A1) fails.
However, Proposition~\ref{prop:jason} can be viewed as a first step to address the question which $D$-finite functions
$f(z)$ define holonomic sequences $(f(n))_{n\geq1}$. The interplay of shift and derivative
in the proof of Proposition~\ref{prop:jason} shows that this problem is not as contrived as it might seem.

\bigskip
{\bf Condition (A3). Growth conditions.} If $f(z)$ is such that $|f(z+1)/f(z)|$ grows faster than 
polynomially as $|z|\to\infty$ then (\ref{leftrightside}) cannot hold identically either. 
We will not use this observation in what follows, since it is well-known that sequences that grow faster than
a power of $n!$ cannot be holonomic.

\bigskip
We now present some sufficient conditions for property~(B).

\bigskip
{\bf Condition (B1). Rational functions.} The class of rational functions has property (B). 
This is trivial but by itself not very useful, since rational functions do not have property (A).
Besides, rational functions clearly define holonomic sequences.

\bigskip
{\bf Condition (B2). Differentiation and Rolle's theorem.}
If $F(z)$ is a non-identically-zero real smooth function
and has infinitely many real positive zeros, then by Rolle's theorem so have $F'(z)$ and all higher derivatives of
$F(z)$. Thus, if for sufficiently large $k$ we get a function $F^{(k)}(z)$ with only finitely many real positive zeros (e.g., a rational function, cf.\ condition~(B1)), then $F(z)$ has only finitely many zeros as well.     
This argument shows that, for example, the class of functions that lie in $\R[x,\log x]$
and are defined in ${[1,{\infty}[}$ has property (B). 

\bigskip
{\bf Condition (B3). Meromorphicity at infinity.} The class of functions $f(z)$ meromorphic at $\infty$ has property (B). Indeed, the function $F(z)=F(z;f,p_0,\dots,p_d)$ is then meromorphic at $\infty$ as well and cannot have infinitely many zeros in $\N$, unless it is identically zero. 
This applies, e.g., to $\sin(1/z)$, $\exp(1/z)$, etc. 

\bigskip
{\bf Condition (B4). Classical results about zeros of analytic functions.}
Many authors have investigated the distribution of zeros of analytic functions.
We just recall the following classical theorem~\cite{PoSz98}. Its gist is that the function $\sin\pi z$
is the ``smallest'' function analytic in the right half-plane that vanishes at $0,1,2,\dots$.
\begin{theorem}[Carlson]\label{thm:carlson}
Let $g(z)$ be a function analytic in $\Re(z)\ge 0$ that satisfies the growth conditions
$g(z)=\O(\exp(\alpha |z|))$ and $g(\pm\im x) = \O(\exp((\pi-\varepsilon)x))$, $x>0$, for some positive constants $\alpha$ and $\varepsilon$.
If $g(z)$ vanishes at $z=0,1,2,\dots$, then $g(z)$ is identically zero. 
\end{theorem}

\begin{corollary}\label{cor:carlson}
  The class of functions $f(z)$ such that $g(z)=f(z+1)$ satisfies the conditions of Theorem~\ref{thm:carlson}
  has property~(B).
\end{corollary}

\bigskip
{\bf Condition (B5). Zeros of real elementary functions -- results of Khovanski\u\i.} In his book
Fewnomials~\cite{Kh91},
Khovanski{\u\i} proves rather general results on the number of simultaneous zeros of multivariate real elementary functions. We will use his results only for univariate functions.
Adopting Khovanski{\u\i}'s definition, we call a univariate function $g(x)$ {\em elementary}, if it can be expressed as a composition of multivariate rational functions and the functions
$\exp(x)$, $\log(x)$, $\sin(x)$, $\cos(x)$, $\tan x$, $\arcsin(x)$, $\arccos(x)$, and $\arctan(x)$.
The domain of definition of $g(x)$ must be such that denominators do not vanish,
arguments of log are positive, arguments of arcsin and arcos lie in $[-1,1]$,
arguments of tan in $[-\pi/2,\pi/2]$, and arguments of sin and cos are bounded. The latter
assumption is crucial to recover the finiteness of the zero set.
Then a special case of one of Khovanski{\u\i}'s main results~\cite[\S1.6]{Kh91} says:
\begin{theorem}[Khovanski{\u\i}]\label{thm:khov}
An elementary function $g(x)$ has only finitely many simple zeros
in its domain of definition.
\end{theorem}

Though apparently only a technical nuisance, the restriction to simple
zeros makes this result somewhat unwieldy for our purposes.
To overcome this difficulty, we impose the additional assumption that $f(x)$ be a quotient
of $D$-finite functions. Then $F(x)$ is also a quotient of $D$-finite functions, and
it is easy to see that the order of the zeros of $F(x)$ is therefore bounded. Thus, by repeated differentiation, we arrive at a function with only simple zeros,
to which Theorem~\ref{thm:khov} can be applied.
\begin{corollary}\label{cor:khov}
The class of elementary functions that are defined on ${[1,{\infty}[}$ and are a quotient of $D$-finite functions has  property~(B).
\end{corollary}
Let us briefly compare this to Carlson's theorem.
Apart from the absence of a growth condition, Khovanski{\u\i}'s result does not require
analyticity in a right half-plane.
Thus, for instance, we need not worry about zeros of the denominator of $f(x)$ outside the positive reals.
On the other hand, Khovanski{\u\i}'s result is about {\em real} functions only, and his restriction to simple zeros
induced us to include the additional $D$-finiteness requirement.


\section{First Results}\label{se:res}

Putting together the results from the preceding section immediately yields several non-holonomicity proofs without
any additional work.

\begin{theorem}\label{thm:log arctan}
  Let $f(z)$ be a function from $\R(z,\log z,\arctan z)$ that is analytic in an open
  set containing the positive reals.
  Then $(f(n))_{n\geq1}$ is a holonomic sequence if and only if $f(z)$ is a rational function.
\end{theorem}
\begin{proof}
  The result is an immediate consequence of conditions~(A1) and (B2) (Corollary~\ref{cor:khov}).
  Indeed, if $f(z)$ depends on $\log z$, then $F(z)=F(z;f,p_0,\dots,p_d)$ has a logarithmic singularity at $z=0$,
  and if it depends on $\arctan z$, then $F(z)$ has a logarithmic singularity at $z=\im$.
\end{proof}
We have excluded the entire functions exp and sin from the preceding theorem, since
condition~(A1) does not work for them in general.
In Section~\ref{se:exp} we will completely
describe the rational functions of $n$ and $\e^n$ that define holonomic sequences.
As long as there is
some log or arctan present, the argument from the proof of Theorem~\ref{thm:log arctan} goes through: For instance,
\[
  a_n = \frac{\log n + \sin \tfrac{1}{n+1}}{\e^n+n^2+1}
\]
does not define a holonomic sequence for this reason.
As for compositions, functions like $\log\log z$ are not allowed in Theorem~\ref{thm:log arctan},
since $D$-finite functions are not closed under composition,
violating our $D$-finiteness assumption that takes care of multiple zeros of elementary functions.
(The sequence $(\log\log n)_{n\geq2}$ can be treated by the asymptotic method~\cite{FlGeSa:05} for showing non-holonomicity.)
Still, $D$-finite functions {\em are} closed under composition with algebraic functions.
For example,
\[
  f(z) = \frac{\e^{z^2}+1}{\e^{z^2}-1}
\]
is a quotient of $D$-finite functions, hence the assumptions of Corollary~\ref{cor:khov} are satisfied.
The complex poles of $f(z)$ are at $z=\pm\sqrt{2k\pi}\e^{\im \pi/4}$, $k\in\Z$, which shows that
condition~(A1) can be applied, establishing the non-holonomicity of the corresponding sequence.
Observe that Carlson's theorem is not (directly) applicable
here, since $f(z)$ is not analytic in a right half-plane. Anyhow, let us now recapitulate what we can 
infer from Corollary~\ref{cor:carlson}.

\begin{theorem}\label{thm:from carlson}
  Let $f(z)$ be a function analytic for $\Re(z)>0$ that does not have an extension that
  is meromorphic at $z=0$. If $f(z+1)$ satisfies the growth conditions
  in Theorem~\ref{thm:carlson}, then the sequence $(f(n))_{n\geq1}$ is not holonomic.
\end{theorem}

This gives another proof that $\log n$ is not holonomic, and also establishes non-holonomicity of $\e^{\sqrt{n}}$, $\e^{1/n}$,
$\exp \exp (1/n)$, and $n^\alpha$ with $\alpha\in\C\setminus\Z$. Note that for the latter three we could have used condition~(A2) instead of (A1)
to establish that $F$ cannot vanish identically. Also, Khovanski\u\i's theorem
can be used instead of Carlson's in all five cases.
Observe that, to our knowledge, this is the first published proof that $\e^{\sqrt{n}}$ and $\exp\exp(1/n)$ are not holonomic.

We have seen by now that our method works on several sequences without any effort.
Take for example the sequence $\e^{1/n}$. Proving its non-holonomicity by asymptotics requires
a somewhat involved asymptotic analysis~\cite{Ge05b} of the generating function $\sum_{n\geq 1}\e^{1/n}z^n$.
But longer proofs can of course have merits, too; this analysis has stimulated comprehensive further asymptotic investigations~\cite{FlGeSa06}
of analytic functions of the form $\sum_{n\geq 1}a^{n^b}z^n$, which are of interest in their own right.
%
%
%
%
%

\section{Factorial and Exponential Sequences}\label{se:exp}

We now present some less obvious applications of the tools we have collected in Section~\ref{se:methods}.
We start with some rather concrete sequences and give a few more general results, most notably
about algebraic sequences, in the following section.
Let us begin with the sequence $n!^\alpha$, $\alpha\in\C\setminus\Z$. It grows too fast to apply Carlson's theorem.
Also, the function $\Gamma(z+1)^\alpha$ is not a quotient of $D$-finite functions,
hence Corollary~\ref{cor:khov} is not applicable. 
We now show that this sequence can be knocked out by Carlson's theorem after a little
bit of massage on the alleged recurrence.

\begin{theorem}
  Let $u_1,\dots,u_s$ be distinct complex numbers,
  and let $\alpha_1,\dots,\alpha_s$ be complex numbers.
  Then the sequence $\Gamma(n-u_1)^{\alpha_1} \dots \Gamma(n-u_s)^{\alpha_s}$
  is holonomic if and only if $\alpha_1,\dots,\alpha_s$ are integers.
\end{theorem}
In particular, powers of hypergeometric sequences with non-integral exponent are
not holonomic, unless the exponent trivially cancels, such as in $(n!^2)^{1/2}$.
This generalizes a known result~\cite[Theorem~1]{Ge04}.

\begin{proof}
  Suppose that the sequence is holonomic.
  Since hypergeometric sequences are holonomic, we may assume w.l.o.g.\ that
  no $\alpha_i$ is an integer.
  The sequence is defined for $n\geq n_0:=\lfloor \sigma \rfloor + 1$, where $\sigma := \max_i \Re(u_i)$.
  Suppose that $\sum_{k=0}^d p_k(z) f(z+k)$ vanishes for $z = n_0,n_0+1,\dots$, where the $p_k(z)$ are polynomials,
  $p_0(z)$ is not identically zero, and
  \[
    f(z) := \Gamma(z-u_1)^{\alpha_1} \dots \Gamma(z-u_s)^{\alpha_s}.
  \]
  By the recurrence of $\Gamma(z)$ and the fact that $\Gamma(z)$
  has no zeros, this implies
  \begin{equation}\label{eq:Ga}
    \sum_{k=0}^d p_k(z) \prod_{i=1}^s \left( (z-u_i)^{\overline{k}} \right)^{\alpha_i} = 0,
     \qquad z = n_0,n_0+1,\dots,
  \end{equation}
  where $z^{\overline{k}}$ denotes the rising factorial
  \[
    z^{\overline{k}} = z(z+1) \dots (z+k-1).
  \]
  Replacing $z$ by $z+n_0$ and applying Theorem~\ref{thm:carlson}, we find
  that the left hand side of~\eqref{eq:Ga} vanishes for $\Re(z)\geq n_0$, hence
  by the identity theorem it vanishes in $\C$ slit along $s$ rays from $u_i$ to $-\infty$.
  We thus have
  \begin{align*}
    0 &= p_0(z) + (z-u_1)^{\alpha_1}\dots(z-u_s)^{\alpha_s} \sum_{k=1}^d p_k(z)
     \prod_{i=1}^s \left( (z-u_i+1)^{\overline{k-1}} \right)^{\alpha_i} \\
    &=: p_0(z) + (z-u_1)^{\alpha_1} G(z).
  \end{align*}
  If we assume w.l.o.g.\ that the $u_i$ are in descending order w.r.t.\ the size of their real parts,
  then $G(z)$ is analytic at $z=u_1$. Since $p_0(z)$ is not identically zero,
  this implies that $(z-u_1)^{\alpha_1}$ is meromorphic at $z=u_1$,
  hence $\alpha_1$ is an integer, a contradiction.
  The converse implication is trivial.
\end{proof}

Another sequence that can be treated in a similar fashion is $n^{\alpha n}$.
So far, its non-holonomicity was only known for rational $\alpha$~\cite{Ge04}.

\begin{proposition}
  For every nonzero complex $\alpha$, the sequence $(n^{\alpha n})_{n\geq1}$ is not holonomic.
\end{proposition}
\begin{proof}
  Let $F(z)$ be as usual. It is easily verified that $(z+1)^{-\alpha(z+1)} F(z+1)$ satisfies the growth
  conditions of Theorem~\ref{thm:carlson}. Appealing to condition~(A1) completes the proof.
\end{proof}

Now let us consider rational functions in $n$ and $\e^n$ with real coefficients. By Khovanski\u\i's theorem
and some commutative algebra, we can precisely characterize the holonomic sequences of this kind,
even if they contain several exponentials.

\begin{theorem}\label{thm:exp}
   Let $\alpha_1,\ldots,\alpha_m$ be positive real numbers such that $\log \alpha_1,\dots,\log \alpha_m$ are
   linearly independent over $\Q$. If $f\in {\mathbb R}(x,\alpha_1^x,\ldots,\alpha_m^x)$
   is such that the sequence $(f(n))_{n\ge 0}$
   is holonomic, then $f(x)\in {\mathbb R}(x)[\alpha_1^{x},\alpha_1^{-x},\ldots, \alpha_m^x,\alpha_m^{-x}]$.
\end{theorem}
\begin{proof}
  We write $\mathbf{y}=(y_1,\dots,y_m)$ for a vector of indeterminates and take some element $g(x,\mathbf{y})=A(x,\mathbf{y})/B(x,\mathbf{y})$
  of $\R(x,\mathbf{y})$ with the property that there is $n_0$ such that
  $(g(n,\alpha_1^n,\dots,\alpha_m^n))_{n\geq n_0}$ is a holonomic sequence.
  W.l.o.g.\ the polynomials $A$ and $B$ are coprime, and $B$ depends on $y_m$.
  By the closure of holonomic sequences under multiplication, we may also assume that $B$ is irreducible.
  Moreover, $B\notin \R[x] \cup \R[x]y_m$, or else there is nothing to show.
  Since this implies that $Ay_m$ and $B$ are coprime, too,
  there are polynomials $L(x,\mathbf{y})$ and $M(x,\mathbf{y})$ such that
  \begin{equation}\label{eq:res}
    L(x,\mathbf{y}) A(x,\mathbf{y})y_m + M(x,\mathbf{y})B(x,\mathbf{y}) = R(x,y_1,\dots,y_{m-1}).
  \end{equation}
  Here, $R$ is the resultant of the two polynomials $Ay_m$ and 
  $B$, viewed as polynomials in the variable $y_m$. Dividing~\eqref{eq:res} by $B$, we may thus assume that the numerator of $g$
  does not depend on $y_m$. Pulling out the coefficient of the highest power of $y_m$, we may at last assume
  that our $g$ is of the form
  \[
    g(x,\mathbf{y}) = \frac{R(x,y_1,\dots,y_{m-1})}{h_0 + \dots + h_{D-1}y_m^{D-1} + y_m^D},
  \]
  where $D>0$, $h_j\in\R(x,y_1,\dots,y_{m-1})$, $h_0\not\equiv 0$, and the denominator is irreducible.
  Suppose now that $g(n,\alpha_1^n,\dots,\alpha_m^n)$ satisfies~\eqref{Prec} for large $n$. By Corollary~\ref{cor:khov}
  (in fact, its obvious extension to ${[n_0,{\infty}[}$), it follows that
  \begin{equation*}\label{eq:poly}
    0 = \sum_{k=0}^d p_k(x) g(x+k,\alpha_1^{x+k},\dots,\alpha_m^{n+k}) = \sum_{k=0}^d p_k(x)g_k(x,\alpha_1^x,\dots,\alpha_m^x)
  \end{equation*}
  for large real $x$. Here, we have defined
  \[
    g_k(x,\mathbf{y}) := g(x+k,\alpha_1^k y_1,\dots,\alpha_m^k y_m).
  \]
  The multiplicative independence of the numbers $\alpha_1,\dots,\alpha_m$ implies the algebraic independence of the functions $\alpha_1^x,\dots,\alpha_m^x$,
  hence
  \[
    \sum_{k=0}^d p_k(x)g_k(x,\mathbf{y}) = 0.
  \]
  Since the numerators of the $g_k$ do not depend on $y_m$, there must be $k>0$ such that the denominators of
  $g_0$ and $g_k$, as polynomials in $\R(x,y_1,\dots,y_{m-1})[y_m]$, have a common root in the algebraic closure
  of $\R(x,y_1,\dots,y_{m-1})$.
  (Essentially, we are once more appealing to condition~(A1).)
  But as these denominators are irreducible, the only way that they can share a root
  is if they coincide, up to a multiplicative ``constant'' from $\R(x,y_1,\dots,y_{m-1})$. The leading coefficients w.r.t.\ $y_m$
  of the denominators of $g_0$ and $g_k$ are $1$ and $\alpha_m^{kD}$, respectively, hence identifying the coefficients of $y_m^0$ yields
  \[
    h_0(x,y_1,\dots,y_{m-1}) = \alpha_m^{-kD} h_0(x,\alpha_1^k y_1,\dots,\alpha_{m-1}^k y_{m-1}).
  \]
  But this equation cannot hold, since comparing the coefficient of any monomial (with non-zero coefficient)
  gives rise to a multiplicative dependence relation between the $\alpha_i$.
\end{proof}

Singer and van der Put~\cite{SivdP97} have shown that the reciprocal $(1/a_n)$ of a holonomic sequence $(a_n)$ is not holonomic
unless $(a_n)$ is a periodic interlacement of hypergeometric sequences.
Thus, in the special case $m=1$ we can close the above proof after showing that the numerator of $g(x,y_1)$ is w.l.o.g.\ independent of $\alpha_1^x$.

To see that the linear independence assumption in Theorem~\ref{thm:exp} is necessary, suppose that
\[
  c_1 \log \alpha_1 + \dots + c_m \log \alpha_m = 0,
\]
where $c_1,\dots,c_m$ are integers, and the $\alpha_i$ are w.l.o.g.\ ordered such that
\[
  c_1,\dots,c_k \geq 0 \qquad \text{and} \qquad c_{k+1},\dots,c_m \leq 0.
\]
Then the sequence $g((n,\alpha_1^n,\dots,\alpha_m^n))_{n\geq 1}$ corresponding to
\[
  g(x,\mathbf{y}) = \frac{1}{1 + y_1^{c_1} \dots y_k^{c_k} - y_{k+1}^{-c_{k+1}} \dots y_m^{-c_m}}
\]
is the constant sequence $(1)_{n\ge1}$.

\section{Miscellaneous Results: Derivatives, Algebraic\\ Functions, and Interlacing}\label{se:misc}

In this section we collect several results, mostly consequences of Khovanski\u\i's theorem.

Take some functions $f(x)$ and $g(x)$
defined for $x>0$ such that $(f(n)_{n\ge1}$ and $(g(n)_{n\ge1}$ are holonomic sequences.
Which operations preserve holonomicity?
Clearly, the sequences arising from $f(x)+g(x)$ and $f(x)g(x)$ are holonomic; we have been using this fact all the time.
Since we assume that $f(x)$ is defined for all positive real numbers, there are some other natural operations
that would not make sense for sequences $(a_n)$ in general.
Integration does not necessarily lead to a function defining a holonomic sequence, as shown by $f(x)=1/x$.
Concerning differentiation, we note the following result.

\begin{proposition}
\label{prop:derivatives}
Let $f(x)$ be an elementary function that is a quotient of $D$-finite functions defined for $x>0$.
If the sequence $(f(n))_{n\ge 1}$ is holonomic, then the sequence $(f'(n))_{n\ge 1}$ is holonomic, too.
\end{proposition}
\begin{proof}
Assume that the sequence $a_n = f(n)$ satisfies~\eqref{Prec}. By Corollary~\ref{cor:khov},
this implies
\[
  f(x+d) = -p_d(x)^{-1} \sum_{k=0}^{d-1}p_k(x)f(x+k), \qquad x\geq 1.
\]
From this relation and Leibniz's rule it follows easily that, for every $\ell\geq d$, the function $f(x+\ell)$ is in the vector space
spanned by $f(x),\dots,f(x+d-1)$ and $f'(x),\dots,f'(x+d-1)$ over $\R(x)$. Since this vector space
has dimension at most $2d$, there must be a linear relation among $f'(x),\dots,f'(x+2d)$.
\end{proof}

Another natural question, which goes slightly beyond closed-form sequences, is whether sequences that satisfy some algebraic equation
can be holonomic. It is clear that periodic interlacements of rational functions are holonomic, but no
other algebraic holonomic sequences are known. For instance, we know from Section~\ref{se:res}
that the sequences $n^\alpha$ with $\alpha\in\Q\setminus\Z$ are not-holonomic, a fact already established
by Gerhold~\cite{Ge04} with a number-theoretic argument.

\begin{theorem}
  Let $f(z)$ be an algebraic function that is analytic in a neighborhood of ${[1,{\infty}[}$. Then
  the sequence $(f(n))_{n\geq 1}$ is holonomic if and only if $f(z)$ is a rational function.
\end{theorem}
\begin{proof}
  Let $F(z)=F(z;f,p_0,\dots,p_d)$ be as usual.  Let $\mathcal{V}$ be the algebraic variety generated by $\{(z,F(z)):z\geq 1\}$.
  By B\'ezout's theorem, the points $(z,F(z))$ stay in one irreducible component of $\mathcal{V}$
  for large real $z$ (otherwise two components would have infinitely many intersections).
  As we suppose that $F(z)$ vanishes for $z\in\N$, this component has infinitely many intersections with the real axis,
  hence it equals the real axis, again by B\'ezout's theorem. But this shows that $F(z)$ vanishes identically.
  Now observe that $f(z)$ has an analytic continuation to a slit plane with finitely many singularities
  and branch points. If all the singularities are poles, then $f(z)$ is rational~\cite[p.~218]{Henrici74}, and
  we are done. Otherwise, we choose a singularity $z_0$, not a pole, such that
  $z_0+1$, $z_0+2,\dots$ are either poles or not singularities of $f(z)$. Then we can appeal
  to condition~(A1), since $F(z)$ is not meromorphic at $z=z_0$.
\end{proof}

A weaker assumption would be that the points $(n,a_n)$ satisfy some algebraic equation.
The result should then be that a holonomic sequence arises only for a periodic interlacement of rational
functions. The problem with the above proof is that the points $(n,a_n)$ may lie in several components
of the variety generated by $\{(n,a_n):n\in\N\}$. We can show as above that each of these
components gives rise to a rational function, but what is missing is an argument
that these rational functions must be interlaced in a periodic way.
The special case where the rational functions are constants is doable, though:

\begin{proposition} 
If a holonomic sequence has only finitely many distinct values,
then it is eventually periodic.
\end{proposition}

\begin{proof}
It is known~\cite[Exercise~VIII.158]{PoSz98b} that this holds if the sequence
satisfies a linear recurrence with {\em constant} coefficients. We show that the result for a sequence
$(a_n)_{n\ge1}$ satisfying a recurrence of the form~\eqref{Prec} follows from this special case.
Let $N$ be the maximum of the degrees of the $p_k$s, and
denote the coefficient of $n^j$ in $p_k(n)$ by $b_{kj}$. Then
\[
  0= \sum_{k=0}^d p_k(n) a_{n+k} = \sum_{j=0}^N L_j(n) n^j
\]
for $n\ge1$, where
\[
  L_j(n) := \sum_{k=0}^d b_{kj}a_{n+k}.
\]
By assumption, the set 
of $N+1$ dimensional vectors
\[
  {\mathcal L}=\{(L_0(n),\ldots,L_N(n)) : n\ge 1\}
\]
is finite.
Since a non-zero polynomial has only finitely many roots, each non-zero vector from $\mathcal{L}$
can occur only for finitely many $n$. In particular, $L_N(n)$ must vanish for large $n$. By definition
of $N$, some $b_{jN}$ is non-zero, which yields a non-trivial constant coefficient recurrence for $(a_n)$.
\end{proof}

If we assume in addition that the values $a_n$ in the preceding proposition are integers, then it follows from the
P{\'o}lya-Carlson theorem~\cite{Ca21b} that $(a_n)$ must be a linear recurrence sequence. This observation is due to
R.P.~Stanley~\cite{St04}.

Finally, we briefly address the question whether non-holonomic sequences can be put together to
produce a holonomic sequence. 

\begin{proposition}
  Suppose that the sequence $a_n$ is a (not necessarily periodic) interlacement of the sequences $\log n$ and $\sqrt{n}$.
  Then $a_n$ is non-holonomic.
\end{proposition}
\begin{proof} 
  Suppose that $a_n$ satisfies a recurrence of order $d$. Let us color the indices $n$ where $\log n$ occurs
  black and the others white. In this coloring of $\N$, there must be some pattern
  of length $d+1$ that occurs infinitely often. Hence one of the $2^{d+1}$ functions
  \[
    F(x) = p_0(x) \genfrac{\{}{\}}{0pt}{}{\sqrt{x}}{\log x} + \dots + p_d(x) \genfrac{\{}{\}}{0pt}{}{\sqrt{x+d}}{\log (x+d)} 
  \]
  has infinitely many zeros. By Corollary~\ref{cor:carlson} or Corollary~\ref{cor:khov}, it vanishes identically, but this can once
  again be refuted by condition~(A1). Hence $a_n$ is not holonomic.
\end{proof}

\section{The Values of the Riemann Zeta Function at Odd Integers}\label{se:zeta}

In this section we consider the sequence $(\zeta(n))_{n\ge2}$, where $\zeta(s)=\sum_{\ell\ge1} s^{-\ell}$ denotes
the Riemann zeta function.
None of our conditions for property~(A) is applicable. In particular, the zeta function
is meromorphic in $\C$, violating condition~(A1). We are thus prompted to look for other arguments.
If $(\zeta(n))_{n\ge2}$ was holonomic, then so would be the subsequence $(\zeta(2n))_{n\ge1}$,
which has the well-known representation 
\[
  \zeta(2n) = (-1)^{n+1}\frac{2^{2n-1}\pi^{2n}B_{2n}}{(2n)!}
\]
in terms of the Bernoulli numbers, whose exponential generating function
\[
  \sum_{n\geq 0} B_{2n} \frac{z^{2n}}{(2n)!} = \frac{z}{\e^z-1}
\]
is the reciprocal of a holonomic function. Holonomic functions are not closed under division, and in fact
Harris and Sibuya~\cite{HaSi85} have shown that the reciprocal $1/f(z)$ of a holonomic function $f(z)$
is holonomic if and only if the logarithmic derivative $f'(z)/f(z)$ is algebraic. From this we infer
that the Bernoulli numbers, hence $(\zeta(2n))_{n\ge1}$, and hence $(\zeta(n))_{n\ge2}$, are not holonomic.

This straightforward argument tells us nothing about the values of the zeta function at \emph{odd} integers. These
numbers have attracted the attention of number theorists for (at least) the last decades, the most well-known
result being Ap\'ery's proof~\cite{Ap79,Poorten79} of the irrationality of $\zeta(3)$. We will now show that the sequence
$(\zeta(2n+1))_{n\ge1}$ is not holonomic. 
Considering the series for $\zeta(2n+1)$ summand by summand, we will deduce a full rank set of linear relations for
the $p_k(n)$.
As an analogy with the theory of transcendental numbers, we mention that
the proof is evocative of the well-known fact that a rapidly converging series can be very useful for´
establishing the transcendence of a number.

\begin{proposition}
  The sequence $(\zeta(2n+1))_{n\ge1}$ is not holonomic.  
\end{proposition}
\begin{proof}
  Suppose that there are polynomials $p_k(n)$ such that
  \[
    \sum_{k=0}^d p_k(n) \zeta(2n+2k+1) = 0, \qquad n\geq 1.
  \]
  By the estimate
  \[
    \sum_{\ell\ge L+1} \ell^{-2n-1} \leq \int_{L}^\infty \frac{\mathrm{d} t}{t^{2n+1}} = \frac{1}{2nL^{2n}} = \O((L-\varepsilon)^{-2n})
  \]
  with arbitrary $\varepsilon>0$, this implies
  \begin{equation}\label{eq:zeta}
    \sum_{k=0}^d p_k(n) = - \sum_{k=0}^d p_k(n) \sum_{\ell\ge2}\ell^{-2n-2k-1} = \O((2-\varepsilon)^{-2n}),
  \end{equation}
  hence the polynomial on the left hand side of~\eqref{eq:zeta} is zero. We iterate this argument as follows: From
  \[
    \sum_{k=0}^d p_k(n)(1 + 2^{-2n-2k-1}) = - \sum_{k=0}^d p_k(n) \sum_{\ell \ge3} \ell^{-2n-2k-1} = \O((3-\varepsilon)^{-2n})
  \]
  we obtain
  \[
    \sum_{k=0}^d \frac{p_k(n)}{2^{2k}} = 0, \qquad n\ge1,
  \]
  and, continuing inductively, we find that
  \[
    \sum_{k=0}^d \frac{p_k(n)}{\ell^{2k}} = 0, \qquad \ell,n\ge1.
  \]
  Since the Vandermonde matrix $(\ell^{-2k})_{1\le \ell \le d+1, 0\le k\le d}$
  is invertible, the polynomials $p_0(n),\dots,p_d(n)$ vanish at all $n\geq1$, hence they vanish identically.
\end{proof}

\section{\bf Comments and Remarks}

%

When trying to establish property~(A), we have to show that a certain function $F$ does not vanish identically.
A possible limit to this endeavor is hinted at by the celebrated results of Richardson~\cite{Richardson68} and Caviness~\cite{Caviness70}.
Recall that they have shown that it is an algorithmically undecidable problem to determine whether a given function from certain classes of elementary functions vanishes identically.
On the one hand, their functions have the absolute value function as a building block, which we do not consider;
on the other hand, they assume a closed form expression, while our $F$ contains arbitrary polynomials.

As for property~(B), we remark that there is an extension of Carlson's classical result due to Malliavin and 
Rubel~\cite{MaRu61,Rubel96}.

Applying Khovanski\u\i's theorem only for quotients of $D$-finite functions severely hampers
the power of the theorem. Maybe the order of zeros of an elementary functions is bounded,
which would make this restriction superfluous.
An observation that might prove useful is that all elementary functions are differentially algebraic, i.e.,
$P(z,f,f',\dots,f^{(k)})=0$ for a nonzero polynomial $P$ (in $k+2$ variables). 
Hopefully this representation can be used to get an upper bound on the order of a zero of $f$. 
However, examples such as $(f')^2+f\cdot f''=0$ show that 
for algebraic differential equations the boundedness cannot be established by the same
straightforward argument as for $D$-finite functions.
Rubel~\cite[Problem 20]{Rubel83}~\cite[Problem 29]{Rubel92} asks which sequences 
$(z_1,z_2,\dots)$ (counting multiplicities) can serve as zero sets of entire differentially algebraic functions 
and remarks that the case $(1,2,3,\dots)$ is still undecided.

Furthermore, we remark that, surprisingly, in the case of integer sequences our restriction to sequences defined via elementary functions is
not that restrictive after all.
Laczkovich and Ruzsa~\cite{LaRu00}
prove that every integer sequence can be represented in the form 
$(f(n))_{n\geq1}$, where $f$ is a so-called naive-elementary function. However, this representation relies on the infinitude
of zeros of the sine function, conflicting with the requirements of Khovanski\u\i's theorem.

We conclude our final remarks with four open problems.
\begin{itemize}
\item[(i)] Besides those described in Theorem~\ref{thm:exp}, which other
elementary (in Khovanski{\u\i}'s sense) functions define holonomic sequences? (Note that $(\sin\alpha n)_{n\geq0}$ clearly is holonomic.)
\item[(ii)] Is a sequence satisfying an algebraic equation holonomic if and only if
it is a periodic interlacement of rational functions?
\item[(iii)] Can our method be pushed forward to accomodate recurrences of more general shapes and multivariate sequences?
\item[(iv)] Let $f(x)$ be a non-rational elementary function such that $(f(n))_{n\ge1}$ is holonomic. If $s(x)$ is
some function such that $(f(s(n)))_{n\ge1}$ is also holonomic, does it follow that $s(x)$ is linear?
\end{itemize}

\bibliographystyle{siam}
\bibliography{../gerhold,../algo}

\end{document}